\documentclass[english]{amsart}
\usepackage[T1]{fontenc}
\usepackage[latin9]{inputenc}
\usepackage{amsthm}
\usepackage{amssymb}

\makeatletter
\numberwithin{equation}{section}
\numberwithin{figure}{section}
\theoremstyle{plain}
\newtheorem{thm}{\protect\theoremname}
\theoremstyle{plain}
\newtheorem{lem}[thm]{\protect\lemmaname}
\newenvironment{lyxcode}
	{\par\begin{list}{}{
		\setlength{\rightmargin}{\leftmargin}
		\setlength{\listparindent}{0pt}
		\raggedright
		\setlength{\itemsep}{0pt}
		\setlength{\parsep}{0pt}
		\normalfont\ttfamily}%
	 \item[]}
	{\end{list}}

\makeatother

\usepackage{babel}
\providecommand{\lemmaname}{Lemma}
\providecommand{\theoremname}{Theorem}

\begin{document}
\title{Complete Monotonicity of Special Functions}
\author{Ruiming Zhang}
\email{ruimingzhang@guet.edu.cn}
\address{School of Mathematics and Computing Sciences\\
Guilin University of Electronic Technology\\
Guilin, Guangxi 541004, P. R. China. }
\thanks{This work is supported by the National Natural Science Foundation
of China grants No. 11771355 and No. 12161022.}
\begin{abstract}
In this work we prove that if an entire function $f(z)$ is of order
strictly less than one and it has only negative zeros, then for each
nonnegative integer $k,m$ the real function $\left(-\frac{1}{x}\right)^{m}\frac{d^{k}}{dx^{k}}\left(x^{k+m}\frac{d^{m}}{dx^{m}}\left(\frac{f'(x)}{f(x)}\right)\right)$
is completely monotonic on $(0,\infty)$. Applications to Askey-Wilson
polynomials, Bessel functions, Ramanujan's entire function, Riemann-xi
function and character Riemann-xi functions are also provided. 
\end{abstract}

\maketitle

\section{\label{sec:Intro} Introduction }

In this work we derive a necessary condition for a genus $0$ entire
function to have only negative zeros, then apply it to an even entire
function of genus $1$ with only real zeros. We apply it to well-known
special functions such as Askey-Wilson polynomials $p_{n}(x;a,b,c,d|q)$,
Bessel functions $J_{\nu}(z),\,K_{iz}(a)$, Ramanujan's entire function
$A_{q}(z)$, $q$-Bessel functions $J_{\nu}^{(2)}(z;q),\,J_{\nu}^{(3)}(z;q)$,
Riemann-xi function $\xi(s)$ and character Riemann-xi function $\xi(s,\chi)$.

Let 

\begin{equation}
f(z)=\sum_{n=0}^{\infty}a_{n}z^{n},\quad a_{0}\cdot a_{n}>0\label{eq:1.1}
\end{equation}
be entire function of order $\rho(f)$ satisfying $0\le\rho(f)<1$
where $\rho(f)$ can be defined by either one of the following formulas,
\begin{equation}
\begin{aligned} & \rho(f)=-\limsup_{n\to\infty}\frac{n\log(n)}{\log|a_{n}|}\\
 & =\limsup_{r\to\infty}\frac{\log\left(\log\left\{ \sup\left\{ \left|f(z)\right|:\left|z\right|\le r\right\} \right\} \right)}{\log r}.
\end{aligned}
\label{eq:1.2}
\end{equation}
Then it is of genus $0$, and it can be represented as
\begin{equation}
\frac{f(z)}{f(0)}=\prod_{n=1}^{\infty}\left(1+\frac{z}{\lambda_{n}}\right),\label{eq:1.3}
\end{equation}
where 
\begin{equation}
\sum_{n=1}^{\infty}\frac{1}{\left|\lambda_{n}\right|}<\infty.\label{eq:1.4}
\end{equation}
Clearly, a genus $1$ even entire function $g$(z) of type

\begin{equation}
g(z)=\sum_{n=0}^{\infty}(-1)^{n}a_{n}z^{2n},\quad a_{0}\cdot a_{n}>0\label{eq:1.5}
\end{equation}
is related to $f(z)$ via the transformation $f(z)=g(i\sqrt{z})$.
Then $g(z)$ has only real zeros $\left\{ \pm z_{n}\right\} _{n\in\mathbb{N}}$
and $g(0)\neq0$, if and only if $f(z)$ has only negative zeros $\left\{ -\lambda_{n}\right\} _{n\in\mathbb{N}}$
with $\lambda_{n}=z_{n}^{2}>0$. 

In this work we prove the following result:
\begin{thm}
\label{thm:1} Let $f(z)$ be an entire function defined in (\ref{eq:1.1})
with $0\le\rho(f)<1$ such that it has only negative roots $\left\{ -\lambda_{n}\right\} _{n\in\mathbb{N}}.$
Then for all nonnegative integers $k,m$ the real function $G_{k}^{(m)}(x)$
defined by
\begin{equation}
G_{k}^{(m)}(x)=\left(-\frac{1}{x}\right)^{m}\frac{d^{k}}{dx^{k}}\left(x^{k+m}\frac{d^{m}}{dx^{m}}\left(\frac{f'(x)}{f(x)}\right)\right),\quad\forall k,m\in\mathbb{N}_{0},\label{eq:1.6}
\end{equation}
are completely monotonic on $(0,\infty)$. The function sequence can
be computed via the following recurrence,
\begin{equation}
G_{0}^{(m)}(x)=(-1)^{m}\frac{d^{m}}{dx^{m}}\left(\frac{f'(x)}{f(x)}\right),\quad G_{k}^{(m)}(x)=(m+k)G_{k-1}^{(m)}(x)+x\frac{d}{dx}G_{k-1}^{(m)}(x),\quad k\in\mathbb{N}.\label{eq:1.7}
\end{equation}
\end{thm}

\section{Proof of Theorem \ref{thm:1}}

\subsection{The heat kernel $\Theta(x)$}
\begin{lem}
\label{lem:2} Assume that $\left\{ \lambda_{n}\right\} _{n\in\mathbb{N}}$
is a sequence of positive numbers such that for certain $0<\alpha<1$
it satisfies
\begin{equation}
0<\sum_{n=1}^{\infty}\frac{1}{\lambda_{n}^{\alpha}}<\infty.\label{eq:2.1}
\end{equation}
 Define 
\begin{equation}
\Theta(x)=\sum_{n=1}^{\infty}e^{-\lambda_{n}x},\quad x>0\label{eq:2.2}
\end{equation}
and $f(z)$ as in (\ref{eq:1.3}). Then $\Theta(x)$ is completely
monotonic on $x\in(0,\infty)$. 

Furthermore, it satisfies 

\begin{equation}
\Theta^{(k)}(x)=\mathcal{O}\left(x^{-\alpha-k}\right),\quad x\to0^{+}\label{eq:2.3}
\end{equation}
and 
\begin{equation}
\Theta^{(k)}(x)=\mathcal{O}\left(e^{-\beta x}\right),\quad x\to+\infty,\label{eq:2.4}
\end{equation}
where $0<\beta<\inf\left\{ \lambda_{n}:n\in\mathbb{N}\right\} $ and
$\alpha$ defined in (\ref{eq:2.1}) with $0<\alpha<1$. 

Moreover, for all $k\in\mathbb{N}_{0},\ x\ge0$,
\begin{equation}
\int_{0}^{\infty}t^{k}e^{-xt}\Theta(t)dt\le\sum_{n=1}^{\infty}\frac{k!}{\lambda_{n}^{k+1}}<\infty\label{eq:2.5}
\end{equation}
and
\begin{equation}
(-1)^{k}\left(\frac{f'(x)}{f(x)}\right)^{(k)}=\sum_{n=1}^{\infty}\frac{k!}{(x+\lambda_{n})^{k+1}}=\int_{0}^{\infty}e^{-xt}t^{k}\Theta(t)dt.\label{eq:2.6}
\end{equation}
\end{lem}

\begin{proof}
Since for any $a>0$
\begin{equation}
\sup_{x>0}x^{a}e^{-x}=\left(\frac{a}{e}\right)^{a},\label{eq:2.7}
\end{equation}
 then for all $t>0,\,k\ge0$ by (\ref{eq:2.7}), (\ref{eq:1.5}) and
(\ref{eq:2.1}), 
\begin{equation}
\begin{aligned} & \sum_{n=1}^{\infty}\lambda_{n}^{k}e^{-\lambda_{n}t}\le\frac{1}{t^{\alpha+k}}\sum_{n=1}^{\infty}\frac{\left(t\lambda_{n}\right)^{k+\alpha}e^{-\lambda_{n}t}}{\lambda_{n}^{\alpha}}\\
 & \le\frac{\sup_{x>0}x^{k+\alpha}e^{-x}}{t^{\alpha+k}}\sum_{n=1}^{\infty}\frac{1}{\lambda_{n}^{\alpha}}=\left(\frac{k+\alpha}{et}\right)^{k+\alpha}\sum_{n=1}^{\infty}\frac{1}{\lambda_{n}^{\alpha}}<\infty,
\end{aligned}
\label{eq:2.8}
\end{equation}
which proves that $\Theta(t)\in C^{\infty}(0,\infty)$ and (\ref{eq:2.3}). 

By (\ref{eq:2.8}) we also have
\begin{equation}
(-1)^{k}\Theta^{(k)}(t)=\sum_{n=1}^{\infty}\lambda_{n}^{k}e^{-\lambda_{n}t}>0,\quad\forall t>0,k\in\mathbb{N}_{0},\label{eq:2.9}
\end{equation}
which proves that $\Theta(t)$ is completely monotonic on $t\in(0,\infty)$
by definition, \cite{SchillingSongVondracek,WidderBook}. 

Let $0<\beta<\inf\left\{ \lambda_{n}:n\in\mathbb{N}\right\} $, then
there exists a positive number $\epsilon$ with $0<\epsilon<1$ such
that 
\begin{equation}
\beta=\epsilon\inf\left\{ \lambda_{n}:n\in\mathbb{N}\right\} \le\epsilon\lambda_{n},\quad\forall n\in\mathbb{N}.\label{eq:2.10}
\end{equation}
For $t\ge1$ by (\ref{eq:2.8}) we have
\begin{equation}
\left|e^{\beta t}\Theta^{(k)}(t)\right|\le\sum_{n=1}^{\infty}\lambda_{n}^{k}e^{-(\lambda_{n}-\beta)t}\le\sum_{n=1}^{\infty}\lambda_{n}^{k}e^{-(1-\epsilon)\lambda_{n}t}\le\sum_{n=1}^{\infty}\lambda_{n}^{k}e^{-(1-\epsilon)\lambda_{n}}<\infty,\label{eq:2.11}
\end{equation}
which establishes (\ref{eq:2.4}).

For all $x,k\ge0$ since
\begin{equation}
0\le\int_{0}^{\infty}t^{k}e^{-xt}\Theta(t)dt=\sum_{n=1}^{\infty}\int_{0}^{\infty}t^{k}e^{-(x+\lambda_{n})t}dt=\sum_{n=1}^{\infty}\frac{k!}{\left(x+\lambda_{n}\right)^{k+1}}\le\sum_{n=1}^{\infty}\frac{k!}{\lambda_{n}^{k+1}}<\infty,\label{eq:2.12}
\end{equation}
then,
\begin{equation}
(-1)^{k}\left(\frac{f'(x)}{f(x)}\right)^{(k)}=\sum_{n=1}^{\infty}\frac{k!}{(x+\lambda_{n})^{k+1}}=\sum_{n=1}^{\infty}\int_{0}^{\infty}t^{k}e^{-(x+\lambda_{n})t}dt=\int_{0}^{\infty}e^{-xt}t^{k}\Theta(t)dt,\label{eq:2.13}
\end{equation}
which proves that for all $k\in\mathbb{N}_{0}$ and $x>0$ we have
\begin{equation}
(-1)^{k}\left(\frac{f'(x)}{f(x)}\right)^{(k)}=\int_{0}^{\infty}e^{-xt}t^{k}\Theta(t)dt\ge0.\label{eq:2.14}
\end{equation}
Equation (\ref{eq:2.14}) establishes that $\frac{f'(x)}{f(x)}$ is
completely monotonic on $(0,\infty)$. 
\end{proof}

\subsection{Proof of the Theorem (\ref{thm:1})}
\begin{proof}
Since $\left\{ \lambda_{n}\right\} _{n\in\mathbb{N}}$ is a positive
sequence, then by Lemma \ref{lem:2} the sequence of functions $\left\{ \Theta_{k}(x)\right\} _{k=0}^{\infty}$
defined by $\Theta_{0}(x)=\Theta(x)$ and
\begin{equation}
\Theta_{k}(x)=(-x)^{k}\Theta^{(k)}(x),\quad\forall x\in(0,\infty),\ \forall k\in\mathbb{N}\label{eq:2.15}
\end{equation}
 are nonnegative functions. Hence, \cite{SchillingSongVondracek,WidderBook}
\begin{equation}
G_{k}^{(m)}(x)=\int_{0}^{\infty}e^{-xt}t^{m}\Theta_{k}(t)dt,\quad\forall k,m\in\mathbb{N}_{0}\label{eq:2.16}
\end{equation}
 are completely monotonic on $x\in(0,\infty)$. By (\ref{eq:2.14})
we get
\begin{equation}
G_{0}^{(m)}(x)=\int_{0}^{\infty}e^{-xt}t^{m}\Theta_{0}(t)dt=\int_{0}^{\infty}e^{-xt}t^{m}\Theta(t)dt=(-1)^{m}\frac{d^{m}}{dx^{m}}\left(\frac{f'(x)}{f(x)}\right).\label{eq:2.17}
\end{equation}
For any $k\in\mathbb{N}$ and $x>0$ apply integration by parts and
(\ref{eq:2.3}) with $0<\alpha<1$, (\ref{eq:2.4}), 
\begin{equation}
\begin{aligned} & G_{k}^{(m)}(x)=\int_{0}^{\infty}e^{-xt}t^{m}\Theta_{k}(t)dt=(-1)^{k}\int_{0}^{\infty}e^{-xt}t^{m+k}\Theta^{(k)}(t)dt\\
 & =(-1)^{k}\int_{0}^{\infty}e^{-xt}t^{m+k}\frac{d}{dt}\Theta^{(k-1)}(t)dt=(-1)^{k-1}\int_{0}^{\infty}\frac{d}{dt}\left(e^{-xt}t^{m+k}\right)\Theta^{(k-1)}(t)dt\\
 & =\dots=\int_{0}^{\infty}\frac{d^{k}}{dt^{k}}\left(e^{-xt}t^{k+m}\right)\Theta(t)dt=x^{-m}\int_{0}^{\infty}\frac{d^{k}}{dy^{k}}\left(e^{-y}y^{k+m}\right)\bigg|_{y=xt}\Theta(t)dt\\
 & =k!\int_{0}^{\infty}e^{-xt}t^{m}\Theta(t)\left(\frac{y^{-m}e^{y}}{k!}\frac{d^{k}}{dy^{k}}\left(e^{-y}y^{k+m}\right)\right)_{y=xt}dt=k!\int_{0}^{\infty}e^{-xt}t^{m}\Theta(t)L_{k}^{(m)}(xt)dt,
\end{aligned}
\label{eq:2.18}
\end{equation}
where we have applied the Rodrigues formula for the Laguerre polynomials
$L_{n}^{(\alpha)}(x)$, \cite{AndrewsAskeyRoy}

\begin{equation}
{\displaystyle L_{n}^{(\alpha)}(x)=\frac{x^{-\alpha}e^{x}}{n!}\frac{d^{n}}{dx^{n}}\left(e^{-x}x^{n+\alpha}\right).}\label{eq:2.19}
\end{equation}
On the other hand, by (\ref{eq:2.18}) and (\ref{eq:2.19}) we have
\begin{equation}
L_{k}^{(m)}(xt)=\frac{x^{-m}e^{xt}}{k!}\frac{d^{k}}{dx^{k}}\left(e^{-xt}x^{k+m}\right)\label{eq:2.20}
\end{equation}
 and
\begin{equation}
\begin{aligned} & G_{k}^{(m)}(x)=k!\int_{0}^{\infty}e^{-xt}t^{m}\Theta(t)L_{k}^{(m)}(xt)dt=x^{-m}\int_{0}^{\infty}t^{m}\Theta(t)\frac{d^{k}}{dx^{k}}\left(e^{-xt}x^{k+m}\right)dt\\
 & =x^{-m}\frac{d^{k}}{dx^{k}}\left(x^{k+m}\int_{0}^{\infty}e^{-xt}t^{m}\Theta(t)dt\right)=\left(-\frac{1}{x}\right)^{m}\frac{d^{k}}{dx^{k}}\left(x^{k+m}\frac{d^{m}}{dx^{m}}\left(\frac{f'(x)}{f(x)}\right)\right),
\end{aligned}
\label{eq:2.21}
\end{equation}
that is, 
\begin{equation}
G_{k}^{(m)}(x)=\left(-\frac{1}{x}\right)^{m}\frac{d^{k}}{dx^{k}}\left(x^{k+m}\frac{d^{m}}{dx^{m}}\left(\frac{f'(x)}{f(x)}\right)\right).\label{eq:2.22}
\end{equation}
For $k\in\mathbb{N}$ it is clear that from (\ref{eq:2.18}) we have
\begin{equation}
\begin{aligned} & G_{k}^{(m)}(x)=(-1)^{k-1}\int_{0}^{\infty}\frac{d}{dt}\left(e^{-xt}t^{m+k}\right)\Theta^{(k-1)}(t)dt\\
 & =(-1)^{k-1}\int_{0}^{\infty}e^{-xt}t^{m+k-1}\left(m+k-xt\right)\Theta^{(k-1)}(t)dt\\
 & =\int_{0}^{\infty}e^{-xt}t^{m}\left(m+k-xt\right)\Theta_{k-1}(t)dt=(m+k)G_{k-1}^{(m)}(x)+x\frac{d}{dx}G_{k-1}^{(m)}(x).
\end{aligned}
\label{eq:2.23}
\end{equation}
This and (\ref{eq:2.17}) gives (\ref{eq:1.7}).
\end{proof}

\section{Apply to the Riemann hypothesis}

\subsection{Orthogonal Polynomials}

Let $\mu(dx)$ be a probability measure on the real line with a support
containing infinitely points such that \cite{Ismail}
\begin{equation}
\int_{\mathbb{R}}\left|x\right|^{k}\mu(dx)<\infty,\quad\forall k\in\mathbb{N}_{0}.\label{eq:3.1}
\end{equation}
 Then there exists a set of polynomials $\left\{ p_{n}(x)\right\} _{n=1}^{\infty}$
satisfying
\begin{equation}
\int_{\mathbb{R}}p_{m}(x)\overline{p_{n}(x)}\mu(dx)=\delta_{m,n},\quad\forall m,n\in\mathbb{N}_{0}.\label{eq:3.2}
\end{equation}
 It is known that for $n\in\mathbb{N}$ all the zeros of $p_{n}(x)$
are real. If the measure is symmetric, $\mu(-dx)=\mu(dx)$, then 
\begin{equation}
p_{2n}(-x)=p_{2n}(x),\quad p_{2n+1}(-x)=-p_{2n+1}(x).\label{eq:3.3}
\end{equation}
Then for all nonnegative integers $k,m$, by Theorem \ref{thm:1},
the following functions 
\begin{equation}
\left(-\frac{1}{x}\right)^{m}\frac{d^{k}}{dx^{k}}\left(x^{k+m}\frac{d^{m+1}}{dx^{m+1}}\log\left(p_{2n}(i\sqrt{x})\right)\right),\quad\left(-\frac{1}{x}\right)^{m}\frac{d^{k}}{dx^{k}}\left(x^{k+m}\frac{d^{m+1}}{dx^{m+1}}\log\left(\frac{p_{2n+1}(i\sqrt{x})}{i\sqrt{x}}\right)\right).\label{eq:3.4}
\end{equation}
are completely monotonic on $(0,\infty)$. Particularly, they are
completely monotonic for $p_{n}(x)$ being the Askey-Wilson polynomials
$p_{n}(x;a,b,c,d)$ or Wilson polynomials $W_{n}(x;a,b,c,d)$, \cite{AndrewsAskeyRoy,KoekoekSwarttouw,Ismail}.

If the support of $\mu(dx)$ is contained in $(0,\infty)$, then for
$n\in\mathbb{N}$ all the roots of $p_{n}(x)$ are in $(0,\infty)$.
Then for all nonnegative integers $k,m$, by Theorem \ref{thm:1}
the functions 
\begin{equation}
\left(-\frac{1}{x}\right)^{m}\frac{d^{k}}{dx^{k}}\left(x^{k+m}\frac{d^{m+1}}{dx^{m+1}}\log p_{n}(-x)\right)\label{eq:3.5}
\end{equation}
are completely monotonic on $(0,\infty)$. For example, if $p_{n}(x)=p_{n}(x-1;a,b,c,d|q),\ n\ge1$
or $p_{n}(x)=Q_{n}(x;\alpha,\beta,N|q)$ the $n$th $q$-Hahn polynomial
for $n=1,\dots,N$, \cite{AndrewsAskeyRoy,Ismail,KoekoekSwarttouw},
then all the functions in (\ref{eq:3.5}) are completely monotonic
on $(0,\infty)$.

\subsection{Bessel functions}

For $\nu>-1$ let $J_{\nu}(z)$ be the Bessel function of the first
kind,\cite{AndrewsAskeyRoy,Gasper1,Gasper2,Ismail}
\begin{equation}
{\displaystyle J_{\nu}(z)=\sum_{m=0}^{\infty}\frac{(-1)^{m}}{m!\,\Gamma(m+\nu+1)}\left(\frac{z}{2}\right)^{2m+\nu},}\label{eq:3.6}
\end{equation}
where $\Gamma(z)$ is the Euler's Gamma function, then the function
$\left(\frac{2}{z}\right)^{\nu}J_{\nu}(z)$ is an even entire function
of order $1$ with only real zeros. Then for all nonnegative integers
$k,m$, by Theorem \ref{thm:1}, the following functions 

\begin{equation}
\left(-\frac{1}{x}\right)^{m}\frac{d^{k}}{dx^{k}}\left(x^{k+m}\frac{d^{m+1}}{dx^{m+1}}\left(\log\frac{I_{\nu}(\sqrt{x})}{(\sqrt{x})^{\nu}}\right)\right)\label{eq:3.7}
\end{equation}
are completely monotonic on $(0,\infty)$ where $I_{\nu}(z)$ is the
modified Bessel function of the first kind,\cite{AndrewsAskeyRoy,Gasper1,Gasper2,Ismail}
\begin{equation}
I_{\nu}(z)=\sum_{m=0}^{\infty}\frac{1}{m!\,\Gamma(m+\nu+1)}\left(\frac{z}{2}\right)^{2m+\nu}.\label{eq:3.8}
\end{equation}
Let $K_{\nu}(z)$ be the modified Bessel functions of second kind
defined by \cite{AndrewsAskeyRoy,Gasper1,Gasper2,Ismail}
\begin{equation}
K_{\nu}(z)=\frac{\pi}{2}\frac{I_{-\nu}(z)-I_{\nu}(z)}{\sin\nu\pi}.\label{eq:3.9}
\end{equation}
 Then for any $a>0$ 
\begin{equation}
K_{iz}(a)=\int_{0}^{\infty}e^{-a\cosh t}\cos(zt)dt,\quad z\in\mathbb{C}\label{eq:3.10}
\end{equation}
is an even entire function of order $1$ that has only zeros. Hence
for all nonnegative integers $k,m$, by Theorem \ref{thm:1}, the
following functions,
\begin{equation}
\left(-\frac{1}{x}\right)^{m}\frac{d^{k}}{dx^{k}}\left(x^{k+m}\frac{d^{m+1}}{dx^{m+1}}\left(\log K_{\sqrt{x}}(a)\right)\right)\label{eq:3.11}
\end{equation}
 are completely monotonic on $(0,\infty)$. 

\subsection{$q$-Transcendental functions}

For $q\in(0,1)$ the $q$-shifted factorial $(a;q)_{n}$ can be defined
by \cite{AndrewsAskeyRoy,Ismail,KoekoekSwarttouw}
\begin{equation}
(a;q)_{\infty}=\prod_{k=0}^{\infty}(1-aq^{k}),\quad(a;q)_{n}=\frac{(a;q)_{\infty}}{(aq^{n};q)_{\infty}},\quad a,n\in\mathbb{C}.\label{eq:3.12}
\end{equation}
 Then the Ramanujan's entire function $A_{q}(z)$, \cite{Ismail}
\begin{equation}
A_{q}(z)=\sum_{n=0}^{\infty}\frac{q^{n^{2}}(-z)^{n}}{(q;q)_{n}},\quad z\in\mathbb{C}\label{eq:3.13}
\end{equation}
 is an entire function of order $0$ with positive zeros. Then for
all nonnegative integers $k,m$, by Theorem \ref{thm:1}, the following
functions,
\begin{equation}
\left(-\frac{1}{x}\right)^{m}\frac{d^{k}}{dx^{k}}\left(x^{k+m}\frac{d^{m+1}}{dx^{m+1}}\left(\log A_{q}(-x)\right)\right)\label{eq:3.14}
\end{equation}
 are completely monotonic on $(0,\infty)$. 

For $\nu>-1$ both 
\begin{equation}
J_{\nu}^{(2)}(z;q)/z^{\nu},\quad J_{\nu}^{(3)}(z;q)/z^{\nu},\quad z\in\mathbb{C}\label{eq:3.15}
\end{equation}
are even entire functions of order $0$ with only real zeros where
the Jackson q-Bessel functions defined by \cite{Jackson1,Jackson2,Jackson3,Jackson4,Jackson5,Koelink}
\begin{align}
J_{\nu}^{(2)}(z;q) & =\frac{(q^{\nu+1};q)_{\infty}}{(q;q)_{\infty}}\left(\frac{z}{2}\right)^{\nu}\sum_{n=0}^{\infty}\frac{q^{n^{2}+n\nu}\left(-z^{2}/4\right)^{n}}{(q,q^{\nu+1};q)_{n}},\label{eq:3.16}\\
J_{\nu}^{(3)}(z;q) & =\frac{(q^{\nu+1};q)_{\infty}}{(q;q)_{\infty}}\left(\frac{z}{2}\right)^{\nu}\sum_{n=0}^{\infty}\frac{q^{n(n+1)/2}\left(-z^{2}/4\right)^{n}}{(q,q^{\nu+1};q)_{n}}.\label{eq:3.17}
\end{align}
Then for all nonnegative integers $k,m$, by Theorem \ref{thm:1},
the following functions,
\begin{equation}
\left(-\frac{1}{x}\right)^{m}\frac{d^{k}}{dx^{k}}\left(x^{k+m}\frac{d^{m+1}}{dx^{m+1}}\left(\log\frac{I_{\nu}^{(2)}(\sqrt{x};q)}{(\sqrt{x})^{\nu}}\right)\right),\quad\left(-\frac{1}{x}\right)^{m}\frac{d^{k}}{dx^{k}}\left(x^{k+m}\frac{d^{m+1}}{dx^{m+1}}\left(\log\frac{I_{\nu}^{(3)}(\sqrt{x};q)}{(\sqrt{x})^{\nu}}\right)\right)\label{eq:3.18}
\end{equation}
 are completely monotonic on $(0,\infty)$ where
\begin{align}
I_{\nu}^{(2)}(z;q) & =\frac{(q^{\nu+1};q)_{\infty}}{(q;q)_{\infty}}\left(\frac{z}{2}\right)^{\nu}\sum_{n=0}^{\infty}\frac{q^{n^{2}+n\nu}\left(z^{2}/4\right)^{n}}{(q,q^{\nu+1};q)_{n}},\label{eq:3.19}\\
I_{\nu}^{(3)}(z;q) & =\frac{(q^{\nu+1};q)_{\infty}}{(q;q)_{\infty}}\left(\frac{z}{2}\right)^{\nu}\sum_{n=0}^{\infty}\frac{q^{n(n+1)/2}\left(z^{2}/4\right)^{n}}{(q,q^{\nu+1};q)_{n}}.\label{eq:3.20}
\end{align}

\subsection{The Riemann $\Xi(s)$ function}

The Riemann $\xi$-function is defined by \cite{AndrewsAskeyRoy,Davenport,Edwards,Gasper1,Gasper2},
\begin{equation}
\xi(s)=\pi^{-s/2}(s-1)\Gamma\left(1+\frac{s}{2}\right)\zeta(s),\quad s=\sigma+it,\ \sigma,t\in\mathbb{R},\label{eq:3.21}
\end{equation}
where $\Gamma(s)$, $\zeta(s)$ are the respective analytic continuations
of
\begin{equation}
\Gamma(s)=\int_{0}^{\infty}e^{-x}x^{s-1}dx,\quad\sigma>0\label{eq:3.22}
\end{equation}
and 
\begin{equation}
\zeta(s)=\sum_{n=1}^{\infty}\frac{1}{n^{s}},\quad\sigma>1.\label{eq:3.23}
\end{equation}
It is well-known that $\xi(s)$ is an order $1$ entire function that
satisfies $\xi(s)=\xi(1-s)$, and it has infinitely many zeros, all
of them are symmetrically distributed inside the vertical strip $\sigma\in(0,1)$,
more than $1/3$ of them on the critical line $\sigma=\frac{1}{2}$.
Then the Riemann Xi function $\Xi(s)=\xi\left(\frac{1}{2}+is\right)$
is an even entire function of order $1$, and it has infinitely many
zeros, all of them are inside the horizontal strip $t\in(-1/2,1/2)$.
The Riemann hypothesis is that all the zeros of $\xi(s)$ are on the
critical line, equivalently, all the zeros of $\Xi(s)$ are real.

Since \cite{Davenport,Edwards,Gasper1,Gasper2}
\begin{equation}
\Xi(s)=\int_{-\infty}^{\infty}\Phi(u)e^{ius}du=2\int_{0}^{\infty}\Phi(u)\cos(us)du,\label{eq:3.24}
\end{equation}
where the even function $\Phi(u)$ satisfies
\begin{equation}
\Phi(u)=\sum_{n=1}^{\infty}\left(4n^{4}\pi^{2}e^{9u/2}-6n^{2}\pi e^{5u/2}\right)e^{-n^{2}\pi e^{2u}}>0\label{eq:3.25}
\end{equation}
on the real line $\mathbb{R}$. Then,
\begin{equation}
\xi\left(\frac{1}{2}+s\right)=\sum_{n=0}^{\infty}a_{n}s^{2n},\quad a_{n}=\frac{2}{(2n)!}\int_{0}^{\infty}\Phi(u)u^{2n}du>0\label{eq:3.26}
\end{equation}
and
\begin{equation}
f(s)=\xi\left(\frac{1}{2}+\sqrt{s}\right)=\sum_{n=0}^{\infty}a_{n}s^{n}\label{eq:3.27}
\end{equation}
defines an entire function of order $\frac{1}{2}$, and $\Xi(s)$
has only real roots if and only if $f(s)$ has only negative roots. 

If the Riemann hypothesis is true, then all the zeros of $f(s)$ are
negative. For all nonnegative integers $m,k$, by Theorem \ref{thm:1},
all the functions,

\begin{equation}
\left(-\frac{1}{x}\right)^{m}\frac{d^{k}}{dx^{k}}\left(x^{k+m}\frac{d^{m+1}}{dx^{m+1}}\left(\log\xi\left(\frac{1}{2}+\sqrt{x}\right)\right)\right)\label{eq:3.28}
\end{equation}
are completely monotonic on $(0,\infty)$. 

\subsection{The character Riemann $\Xi(s,\chi)$ function}

Given a primitive Dirichlet character $\chi(n)$ modulo $q$, let
\cite{AndrewsAskeyRoy,Davenport,Edwards}
\begin{equation}
\xi(s,\chi)=\left(\frac{q}{\pi}\right)^{(s+\kappa)/2}\Gamma\left(\frac{s+\kappa}{2}\right)L(s,\chi),\label{eq:3.29}
\end{equation}
where $\kappa$ is the parity of $\chi$ and $L(s,\chi)$ is the analytic
continuation of 
\begin{equation}
L\left(s,\chi\right)=\sum_{n=1}^{\infty}\frac{\chi(n)}{n^{s}},\quad\sigma>1.\label{eq:3.30}
\end{equation}
Then $\xi(s,\chi)$ is an entire function of order $1$ such that
\cite{Davenport}
\begin{equation}
\xi(s,\chi)=\epsilon(\chi)\xi(1-s,\overline{\chi}),\label{eq:3.31}
\end{equation}
where 
\begin{equation}
\epsilon(\chi)=\frac{\tau(\chi)}{i^{\kappa}\sqrt{q}},\quad\tau(\chi)=\sum_{n=1}^{q}\chi(n)\exp\left(\frac{2\pi in}{q}\right).\label{eq:3.32}
\end{equation}
Let 
\begin{equation}
G(s,\chi)=\xi\left(s,\chi\right)\cdot\xi\left(s,\overline{\chi}\right),\label{eq:3.33}
\end{equation}
then
\begin{equation}
G(s,\chi)=\epsilon\left(\chi\right)\cdot\epsilon\left(\overline{\chi}\right)G(1-s,\chi).\label{eq:3.34}
\end{equation}
Since \cite{Davenport}
\begin{equation}
\tau\left(\overline{\chi}\right)=\overline{\tau(\chi)},\quad\left|\tau(\chi)\right|=\sqrt{q},\label{eq:3.35}
\end{equation}
then
\begin{equation}
G(s,\chi)=G(1-s,\chi).\label{eq:3.36}
\end{equation}
Since the entire function $\xi\left(\frac{1}{2}+is,\chi\right)$ has
an integral representation, \cite{Davenport}
\begin{equation}
\xi\left(\frac{1}{2}+is,\chi\right)=\int_{-\infty}^{\infty}e^{isy}\varphi\left(y,\chi\right)dy,\label{eq:3.37}
\end{equation}
where
\begin{equation}
\varphi(y,\chi)=2\sum_{n=1}^{\infty}n^{\kappa}\chi(n)\exp\left(-\frac{n^{2}\pi}{q}e^{2y}+\left(\kappa+\frac{1}{2}\right)y\right),\label{eq:3.38}
\end{equation}
then
\begin{equation}
\xi\left(\frac{1}{2}+is,\chi\right)=\sum_{n=0}^{\infty}i^{n}a_{n}(\chi)s^{n},\quad a_{n}(\chi)=\int_{-\infty}^{\infty}y^{n}\varphi\left(y,\chi\right)dy.\label{eq:3.39}
\end{equation}
It is known that the fast decreasing smooth function $\varphi(y,\chi)$
satisfies the functional equation \cite{Davenport}
\begin{equation}
\varphi(y,\chi)=\frac{i^{\kappa}\sqrt{q}}{\tau\left(\overline{\chi}\right)}\varphi(-y;\overline{\chi}),\quad y\in\mathbb{R},\label{eq:3.40}
\end{equation}
then for all $n\in\mathbb{N}_{0}$,
\begin{equation}
\begin{aligned} & a_{n}(\overline{\chi})=\int_{-\infty}^{\infty}y^{n}\varphi\left(y,\overline{\chi}\right)dy=(-1)^{n}\int_{-\infty}^{\infty}y^{n}\varphi\left(-y,\overline{\chi}\right)dy\\
 & =\frac{(-1)^{n}\tau\left(\overline{\chi}\right)}{i^{\kappa}\sqrt{q}}\int_{-\infty}^{\infty}y^{n}\varphi\left(y,\chi\right)dy=\frac{(-1)^{n}\tau\left(\overline{\chi}\right)}{i^{\kappa}\sqrt{q}}a_{n}(\chi).
\end{aligned}
\label{eq:3.41}
\end{equation}
Let
\begin{equation}
f(s,\chi)=s^{-2\mu}G\left(\frac{1}{2}+is,\chi\right)=s^{-2\mu}\xi\left(\frac{1}{2}+is,\chi\right)\cdot\xi\left(\frac{1}{2}+is,\overline{\chi}\right),\label{eq:3.42}
\end{equation}
where $\mu\in\mathbb{N}_{0}$ is the least nonnegative integer such
that $a_{\mu}(\chi)\neq0$. Then the even entire function $f(s,\chi)$
has the series expansion
\begin{equation}
f(s,\chi)=f(-s,\chi)=\sum_{n=0}^{\infty}(-1)^{n}b_{n}(\chi)s^{2n},\label{eq:3.43}
\end{equation}
where $b_{0}(\chi)=\frac{(-1)^{\mu}\tau\left(\overline{\chi}\right)}{i^{\kappa}\sqrt{q}}a_{\mu}^{2}(\chi)$
and for $n\in\mathbb{N}$,
\begin{equation}
\begin{aligned} & b_{n}(\chi)=\sum_{j=0}^{2n}a_{j+\mu}(\chi)a_{2n-j+\mu}(\overline{\chi})\\
 & =\frac{(-1)^{\mu}\tau\left(\overline{\chi}\right)}{i^{\kappa}\sqrt{q}}\sum_{j=0}^{2n}(-1)^{j}a_{j+\mu}(\chi)a_{2n-j+\mu}(\chi).
\end{aligned}
\label{eq:3.44}
\end{equation}
It is well-known that $\xi(s,\chi)$ has infinitely many zeros, all
of them are in the horizontal strip $t\in(-1/2,1/2)$, \cite{Davenport}.
Then $f(s,\chi)$ is an order $1$ even entire function with infinitely
many zeros, all of them are in the horizontal strip $t\in(-1/2,1/2)$.
Clearly, all the zeros of $\xi(s,\chi)$ on the critical line $\sigma=\frac{1}{2}$
if and only if all the zeros of $f(s,\chi)$ are real. Therefore,
the generalized Riemann hypothesis for $L\left(s,\chi\right)$ is
equivalent to that all the zeros of $f(s,\chi)$ are real. 

For all nonnegative integers $m,k$, assuming the generalized Riemann
hypothesis associated with $\chi$, we apply Theorem \ref{thm:1}
to $f(s,\chi)$ to have all the real functions,

\begin{equation}
\left(-\frac{1}{x}\right)^{m}\frac{d^{k}}{dx^{k}}\left(x^{k+m}\frac{d^{m+1}}{dx^{m+1}}\left(\log f\left(\sqrt{x},\chi\right)\right)\right)\label{eq:3.45}
\end{equation}
are completely monotonic on $(0,\infty)$. 

Formulas (\ref{eq:3.28}) and (\ref{eq:3.45}) could be used to disprove
the corresponding Riemann hypothesis.
\begin{lyxcode}
\end{lyxcode}


\begin{thebibliography}{10}
\bibitem{Ahlfors}L. Ahlfors, \emph{Complex Analysis}, 3rd edition,
McGraw-Hill Education, 1979.

\bibitem{AndrewsAskeyRoy}G. Andrews, R. Askey and R. Roy, \emph{Special
Functions}, 1st edition, Cambridge University Press.

\bibitem{Boas}R. P. Boas, \emph{Entire Functions}, 1st edition, Academic
Press, 1954.

\bibitem{Davenport}H. Davenport, \emph{Multiplicative Number Theory},
Springer-Verlag, New York, 1980.

\bibitem{Edwards}H. M. Edwards, \emph{Riemann's Zeta Function}, Dover
Publications, 2001.

\bibitem{Gasper1}George Gasper, Using sums of squares to prove that
certain entire functions have only real zeros, dedicated to the memory
of Ralph P. Boas, Jr. (1912--1992), published in \emph{Fourier Analysis:
Analytic and Geometric Aspects}, W.O. Bray, P.S. Milojevic and C.V.
Stanojevic, eds., Marcel Dekker, 1994, pp. 171--186.

\bibitem{Gasper2} G. Gasper, Using integrals of squares to prove
that certain real-valued special functions to prove that the P\'{o}lya
$\Xi^{*}\left(z\right)$ function, the functions $K_{iz}\left(a\right),\ a>0$
and some other entire functions have only real zeros, dedicated to
Dan Waterman on the occasion of his 80th birthday.

\bibitem{Ismail}M. E. H. Ismail, \emph{Classical and Quantum Orthogonal
Polynomials in One Variable, }Cambridge University Press, Cambridge,
2005.

\bibitem{Jackson1}F. H. Jackson, I.---On generalized functions of
Legendre and Bessel, Transactions of the Royal Society of Edinburgh,
41 (1): 1--28, 1906.

\bibitem{Jackson2}F. H. Jackson, VI.---Theorems relating to a generalization
of the Bessel function, Transactions of the Royal Society of Edinburgh,
41 (1): 105--118, 1906.

\bibitem{Jackson3} F. H. Jackson, XVII.---Theorems relating to a
generalization of Bessel's function, Transactions of the Royal Society
of Edinburgh, 41 (2): 399--408, 1906.

\bibitem{Jackson4} F. H. Jackson, The Application of Basic Numbers
to Bessel's and Legendre's Functions, Proceedings of the London Mathematical
Society, 2, 2 (1): 192--220, 1905. 

\bibitem{Jackson5}F. H. Jackson, The Application of Basic Numbers
to Bessel's and Legendre's Functions (Second paper), Proceedings of
the London Mathematical Society, 2, 3 (1): 1--23, 1905. 

\bibitem{Koelink}H. T. Koelink, Hansen-Lommel Orthogonality Relations
for Jackson's q-Bessel Functions, Journal of Mathematical Analysis
and Applications, 175 (2): 425--437, 1993.

\bibitem{KoekoekSwarttouw}R. Koekoek and R. F. Swarttouw, \emph{The
Askey-scheme of hypergeometric orthogonal polynomials and its q-analogue},
Delft University of Technology, Faculty of Information Technology
and Systems, Department of Technical Mathematics and Informatics,Report
no. 98-17, 1998. https://fa.ewi.tudelft.nl/\textasciitilde koekoek/askey/.

\bibitem{SchillingSongVondracek}R. Schilling, R. Song and Z. Vondracek,
\emph{Bernstein Functions: Theory and Applications}, 1st edition,
de Gruyter, 2010.

\bibitem{WidderBook}D. V. Widder, \emph{The Laplace Transform}, Princeton
Mathematical Series, v. 6, Princeton University Press, 1941.
\end{thebibliography}
\end{document}